\documentclass{amsart} 
\usepackage{amssymb, amsmath}

   \def\CC{\mathbb C}
   
   \def\NN{\mathbb N}
   \def\ZZ{\mathbb Z}
   \def\RR{\mathbb R}
   \def\cS{\mathcal S}
   \def\om{\omega}
   
   \def\la{\lambda}
   \def\al{\alpha}
   
   \def\ds{\displaystyle}
    
   \def\Re{\mathop{\rm Re}\nolimits}
   \newcommand {\norm}[1] {\| #1 \|}
   \newcommand {\lrnorm}[1]{\left\| #1 \right\|}
   \newcommand {\bignorm}[1]{\bigl\| #1 \bigr\|}
   \newcommand {\Bignorm}[1]{\Bigl\| #1 \Bigr\|}
   \newcommand {\biggnorm}[1]{\biggl\| #1 \biggr\|}
   \newcommand {\sfrac}[2] { {{}^{#1}\!\!/\!{}_{#2}}} 
   \newcommand {\onehalf}{\sfrac12}

  \theoremstyle{plain}
  \newtheorem{lemma}{Lemma}[section]
  \newtheorem{theorem}[lemma]{Theorem}
  \newtheorem{proposition}[lemma]{Proposition}
  
  \theoremstyle{definition}
  \newtheorem{definition}[lemma]{Definition}
  \newtheorem{remark}[lemma]{Remark}
  \newtheorem{example}[lemma]{Example}
  \allowdisplaybreaks

  \renewcommand{\labelenumi} {(\alph{enumi})}    
  \renewcommand{\labelenumii}{(\roman{enumii})}
  \renewcommand{\theenumi} {(\alph{enumi})}      
  \renewcommand{\theenumii}{(\roman{enumii})}

\begin{document}

\title[On diagonal Volterra equations]{Admissibility and Controllability of diagonal Volterra
  equations with scalar inputs}

\author[Haak]{Bernhard H.~Haak}
\address{Institut de Math\'ematiques de Bordeaux, 
                    Universit\'e Bordeaux 1, 
                    351 cours de la Lib\'eration, %
                    F-33405 Talence cedex, %
                    France} 
\email{Bernhard.Haak@math.u-bordeaux1.fr}

\author[Jacob]{Birgit Jacob}
\address{Institute for Mathematics, %
                    University of Paderborn, %
                    Warburger Str. 100, %
                    D--33098 Paderborn, %
                    Germany}
\email{jacob@math.uni-paderborn.de}

\author[Partington]{Jonathan R.~Partington}
\address   {School of Mathematics, %
                    University of Leeds, %
                    Leeds LS2 9JT, %
                    U.K.}
\email{J.R.Partington@leeds.ac.uk}

\author[Pott]{Sandra Pott}
\address{Department of Mathematics, %
                    University of Glasgow, %
                    Glasgow, G12 8QW, %
                    U.K. \\and \\
Institute for Mathematics, %
                    University of Paderborn, %
                    Warburger Str. 100, %
                    D--33098 Paderborn, %
                    Germany}
\email{sandrapo@math.uni-paderborn.de}

\let\cal\mathcal
\thanks{The authors are partially supported by the Royal Society's
  {\em International Joint Project\/} scheme and from the German
  Research Foundation (DFG) which they gratefully acknowledge}

\begin{abstract}
  This article studies Volterra evolution equations from the point of
  view of control theory, in the case that the generator of the
  underlying semigroup has a Riesz basis of eigenvectors. Conditions
  for admissibility of the system's control operator are given in
  terms of the Carleson embedding properties of certain discrete
  measures. Moreover, exact and null controllability are expressed in
  terms of a new interpolation question for analytic functions,
  providing a generalization of results known to hold for the standard
  Cauchy problem. The results are illustrated by examples involving heat 
  conduction with memory.
\end{abstract}
\def\sep{, }
\keywords{Control of Volterra equations \sep Admissibility \sep 
Controllability \sep Carleson measures}

  \sloppy
  \allowdisplaybreaks

  \renewcommand{\labelenumi} {(\alph{enumi})}    
  \renewcommand{\labelenumii}{(\roman{enumii})}
  \renewcommand{\theenumi} {(\alph{enumi})}      
  \renewcommand{\theenumii}{(\roman{enumii})}

\maketitle

\section{Introduction}
Consider the evolution equation
\begin{equation}\label{eqn:controlled}
 x(t) = x_0+\int_0^t a(t-s) A x(s)\,ds +\int_0^t Bu(s)\, ds, \qquad t\ge 0.
\end{equation}
Here we assume that $A$ generates a $C_0$-semigroup $(T(t))_{t\ge 0}$
on a Hilbert space $X$, $a \in L^1_{loc}(0,\infty)$ is real-valued and
of at most exponential growth, and the control operator $B\in {\cal
  L}(U, D(A^*)^*)$, where $U$ is another Hilbert space. It is further
assumed that the uncontrolled system
\begin{equation}\label{eqn:uncontrolled}
 x(t) = x_0+\int_0^t a(t-s) A x(s)\,ds  \qquad t\ge 0
\end{equation}
is well-posed, which is equivalent to the existence of a unique family
of bounded linear operators $(S(t))_{t\ge 0}$ on $X$, such that
\begin{enumerate}
\item $S(0)=I$ and  $(S(t))_{t\ge 0}$ is strongly continuous on $\mathbb R_+$.
\item $S(t)$ commutes with $A$, which means $S(t)(D(A))\subset D(A)$ for all
  $t\ge 0$, and $AS(t)x=S(t)Ax$ for all $x\in D(A)$  and $t\ge 0$.
\item For all $x\in D(A)$ and all $t\ge 0$ the resolvent equations hold:
\begin{equation}\label{eqn:resolvent}
S(t)x = x +\int_0^t a(t-s)AS(s)x\,ds.
\end{equation}
\end{enumerate}
The family of bounded linear operators $(S(t))_{t\ge 0}$ is called the
{\em resolvent} or {\em solution family} for (\ref{eqn:uncontrolled}).
We refer to the monograph by Pr\"uss \cite{pr} for more about
resolvents.  In particular, if we assume further that the resolvent
$(S(t))_{t\ge 0}$ is exponentially bounded, say $\norm{ S(t) } \le
Me^{\omega t}$ for $t \ge 0$, then the Laplace transform of
$S(\cdot)x_0$ is well-defined and satisfies
\[
\widehat S(\la)x_0= \frac{1}{\la} (I-\widehat a(\la) A)^{-1} x_0 
   \qquad (\Re \la>\omega)
\]
(here the hat denotes Laplace transform). The assumption of an
exponential growth of the resolvent is indeed a restriction of
generality: in contrast with semigroups or cosine families, resolvents
may grow super-exponentially in time even if the kernel $a$ is
integrable and of class $C^\infty$ (see \cite{DeschPruess} for more
details).\\

Notice that by adding $\om\cdot a\ast x$ on both sides of
equation~(\ref{eqn:uncontrolled}) we obtain an equation of the same
form where $x$ is replaced by $v=x+\om\cdot a\ast x$, $A$ is replaced
by $A+\om$, and $a$ by the solution $r$ of $r+\om\cdot a\ast r =
a$. Indeed,
\begin{eqnarray*}
 v &=& x_0 + a\ast (A+\om)x = x_0 + [r+\om\cdot a\ast r]\ast (A+\om)x \\
   &=& x_0+ r\ast(A+\om)x + \om\cdot a\ast r\ast(A+\om)x \\
   &=& x_0 + r\ast (A+\om)(x+\om\cdot a\ast x) = x_0+ r\ast(A+\om)v. 
\end{eqnarray*}
This transformation shows that without loss of generality we may
assume $A$ to generate a uniformly exponentially stable semigroup.  We
notice that $1/\widehat r(\la) = [1/\widehat a(\la)] + \om$ in this
case.

\begin{example}\label{ex:important-kernels}
\begin{enumerate}
\item Consider the standard kernel $a(t) =
  \frac{t^{\beta-1}}{\Gamma(\beta)}$ for $\beta \in (0,2)$ given in
  \cite[Example 2.1]{pr}.  We have $\widehat
  a(\lambda)=\lambda^{-\beta}$. In our main result in
  Theorem~\ref{thm:admissibility-result} we consider a class of
  kernels that admit upper and lower estimates against this standard
  kernel.
\item Another important class of kernels is that given by
  \cite[Example 2.2]{pr}:
  \begin{equation}
  a(t)= \int_0^\infty e^{-st} \, d\alpha(s), \label{eq:class2}
  \end{equation}
  or
  \[
  \widehat a(\lambda)= \int_0^\infty \frac{1}{\lambda+s} \, d\alpha(s),
  \]
  where $\alpha$ is a non-decreasing function on $[0,\infty)$ such that
  \[
  \int_1^\infty d\alpha(s)/s < \infty. 
  \]
  In Section~\label{item:1}~\ref{sec:controllability} we give a result
  on controllability of special cases of such kernels.
  \item Let $a(t) = \int_0^\infty \frac{t^{\rho-1}}{\Gamma(\rho)}
    \,d\rho$ as considered in \cite[Example 2.3]{pr}. We then have 
     $\widehat a(\la)=1/\log(\la)$. In Theorem~\ref{thm:admissibility-result}
     we obtain a sufficient criterion for admissibility in this case.
  \end{enumerate}
\end{example}

The mild solution of (\ref{eqn:controlled}) is formally 
given by the variation of constants formula
\[
  x(t)=S(t)x_0 +(S*Bu)(t),\qquad t\ge 0,
\]
which is actually the classical solution if $B\in {\cal L}(U,X)$,
$x_0\in D(A)$ and $u$ sufficiently smooth. In general however, $B$ is
not a bounded operator from $U$ to $X$ and so an additional assumption
on $B$ will be needed to ensure $x(t)\in X$ for every $x_0\in X$ and
every $u\in L^2(0,\infty;\,U)$.\\

In Section~\ref{sec:2} we introduce the idea of admissibility for a
control operator, and explain some of its properties. We also present
the more familiar theme of controllability, which will also be studied
in this paper.  The main results of the paper are contained in
Section~\ref{sec:3}, where we specialise to diagonal systems, and
derive conditions for admissibility and controllability of such
systems, presented in terms of Carleson embedding and interpolation
properties.  Finally, in Section \ref{sec:examples}, we illustrate the
ideas of this paper with examples involving heat conduction with
memory.

\section{Admissibility and Controllability}
\label{sec:2}

Since the resolvent for (\ref{eqn:uncontrolled}) commutes with the
operator $A$, it can be easily seen that the resolvent operator
$(S(t))_{t\ge 0}$ can be restricted/extended to a resolvent operator
on $D(A)$/$D(A^*)^*$. We denote the restriction/extension again by
$(S(t))_{t\ge 0}$.  Similarly, the operator $A$ can be
extended/restricted to a generator of a $C_0$-semigroup on
$D(A)$/$D(A^*)^*$, again denoted by $A$.

\begin{definition}\label{def:admissibility}
Let $B\in{\cal L}(U,D(A^*)^*)$. Then $B$ is called {\em admissible
for $(S(t))_{t\ge 0}$} if there exists a constant $M>0$ such that
\begin{equation}  \label{eq:admiss-a}
  \norm{ (S*Bu)(t) }_{X} \le M \norm{ u }_{L^2(0,\infty;\,U)},
      \qquad u\in L^2(0,\infty;\,U),\qquad t\ge 0.
\end{equation}
\end{definition}

\begin{remark}\label{rem:admiss-equiv}
$B\in{\cal L}(U,D(A^*)^*)$ is admissible if and only if
there exists a constant $M>0$ such that
\begin{equation}  \label{eq:admiss-b}
 \lrnorm{  \int_0^\infty S(t)Bu(t) \, dt }_X \le M \norm{ u }_{L^2(0,\infty;\,U)}.
\end{equation}
for all $u \in L^2(0,\infty;\, U)$ with compact support.
\end{remark}
\begin{proof}
  Let $B$ be admissible and assume that $u \in L^2(0,\infty; \,U)$ has
  compact support, say $[a,b]$. Define $\widetilde u(s) := u(b-s)$
  when $b{-}s\in [a,b]$ and zero otherwise. Then,
\[
    \Bignorm{ \int_0^\infty S(s)B u(s)\,ds }
=   \bignorm{ (S \ast B\widetilde u) (b)  }
\le M  \bignorm{ \widetilde u }_{L^2(0,\infty;\,U)}
=   M \bignorm{ u }_{L^2(0,\infty;\,U)} .
\]
Conversely, let $t\ge0$ and $u\in L^2(0,\infty;\, U)$. Define $v_t(s)
= u(t{-}s)$ for $s\in [0,t]$ and zero otherwise. Then
(\ref{eq:admiss-b}) implies
\[
\norm{  (S\ast B u)(t)  } = \Bignorm{  \int_0^\infty S(s)B v_t(s)\,ds }
\le M \norm{  v_t  }_{L^2(0,t;\, U)} \le M \norm{  u  }_{L^2(0,\infty;\, U)},
\]
whence $B$ is admissible. 
\end{proof}

Admissibility of the operator $B$ guarantees that the operator 
\[{\cal B}_{\infty}:\{u\in L^2(0,\infty;\,U)\mid 
u \mbox{ has compact support}\} \rightarrow X,\]
 given by,
\begin{equation}\label{binfty} 
   {\cal B}_{\infty}u := \int_0^\infty S(s)B u(s)\,ds,
\end{equation}
possesses a unique extension to a linear, bounded operator from
$L^2(0,\infty;\,U)$ to $X$. We denote this extension again by ${\cal
  B}_{\infty}$. If the solution family is exponentially stable, then
formula (\ref{binfty}) holds for
every $u\in L^2(0,\infty;\,U)$.\\

There is also the notion of admissibility of an observation operator
$C \in {\cal L}(D(A),Y)$, where $Y$ is another Hilbert space,
guaranteeing that the output $y$, where
\[
y(t) = C x(t), \qquad t \ge 0,
\]
lies in $L^2$. For infinite-time admissibility, the following is the
most natural definition.

\begin{definition}
  The operator $C$ is called an {\em admissible observation operator}
  for the uncontrolled system (\ref{eqn:uncontrolled}), if there
  exists a constant $M>0$ such that
\[
    \norm{ y(\cdot) }_{L^2(0,\infty;Y)}
  = \norm{ CS(\cdot)x_0 }_{L^2(0,\infty;Y)} 
\le M \norm{ x_0 },\qquad x_0\in D(A).
\]
The operator $C$ is called a {\em finite-time admissible observation
operator} for (\ref{eqn:uncontrolled}), if there exist constants
$M>0$ and $\omega\in\mathbb R$ such that
\[
   \norm{ y(\cdot) }_{}
 = \norm{ CS(\cdot)x_0 }_{L^2(0,t;Y)} 
\le Me^{\omega t} \norm{ x_0 },\qquad x_0\in D(A),t>0.
\]
\end{definition}

Notice that the dual operator ${\cal B}^{\;*}_\infty$ is given by 
$x^* \mapsto B^* S(\cdot)^* x^*$. Therefore, there is a natural duality
between admissibility of control operators and admissibility of observation
operators, that is, $B \in {\cal L}(U,D(A^*)^*)$ is an admissible
control operator if and only if $B^* \in {\cal L}(D(A^*),U^*)$ is an
admissible observation operator. This is explained in detail in 
\cite[Section~4]{jp04}. For more on admissibility for the
Cauchy problem (i.e., $a \equiv 1$), we refer to the survey \cite{jpsurvey}.

We shall also be interested in obtaining conditions for exact
controllability of the system (\ref{eqn:controlled}). Accordingly, we
make the following definitions.

\begin{definition}
  The system (\ref{eqn:controlled}) is said to be {\em exactly
    controllable}, if every state can be achieved by a suitable
  control, i.e., if $R({\cal B}_{\infty}) \supseteq X$.

It is said to be {\em null-controllable in time\/} $\tau>0$ if 
$R({\cal  B}_{\infty}) \supseteq R(S(\tau))$.
\end{definition}

For a recent discussion of these properties in the context of the
Cauchy problem, we refer to \cite{jp06}.

\section{Admissible and Controllable Diagonal Systems}
\label{sec:3}
{}From now on we assume that $A$ is the infinitesimal generator of an
exponentially stable $C_0$-semigroup $(T(t))_{t\ge 0}$ on $H$ with a
sequence of normalised eigenvectors $\{\phi_n\}_{n\in \mathbb N}$
forming a Riesz basis for $H$, with associated eigenvalues
$\{\lambda_n\}_{n\in\mathbb N}$, that is,
\[
   A\phi_n= \lambda_n\phi_n, \qquad n\in\mathbb N.
\]
Let $\cS(\theta)$ be the open sector of angle $2\theta$ symmetric
about the positive real axis.
Recall that the condition $-\la_n \in \cS(\sfrac\pi2)$ for all $n\in
\NN$ is necessary for $A$ to generate a bounded semigroup and that
$-\lambda_n \in \cS(\theta)$ with $\theta < \sfrac\pi2$ is equivalent
to $A$ generating a bounded {\em analytic} semigroup.

\medskip

Since $(T(t))_{t\ge 0}$ is assumed to be exponentially stable we have
$\sup_{n\in\mathbb N}\mbox{Re}\,\lambda_n<0$. Let $\psi_n$ be an
eigenvector of $A^*$ corresponding to the eigenvalue
$\overline{\lambda_n}$. Without loss of generality we can assume that
$\langle \phi_n,\psi_n\rangle =1$. Then the sequence
$\{\psi_n\}_{n\in\mathbb N}$ forms a Riesz basis of $H$ and every
$x\in H$ can be written as 
\[ 
   x=\sum_{n\in\mathbb N} \langle x,\psi_n\rangle
   \phi_n=\sum_{n\in\mathbb N} \langle x,\phi_n\rangle \psi_n. 
\]
Note that the Volterra system is also diagonal, in the 
sense that  there exist functions  $c_n$ such that
$S(t)\phi_n=c_n(t) \phi_n$; indeed
\[
   \widehat S(\la) \phi_n = \sigma(\la, -\lambda_n) \phi_n,
\]
where
\[
   \sigma(\lambda, \mu) = \frac{1}{\la(1+\mu \widehat a(\lambda))},
   \qquad \Re  \mu,\Re \lambda>0,
\]
so the Laplace transform of $c_n$ is $\sigma(\cdot,-\lambda_n)$.\\

For example, if $\widehat a(\la)=\xi/(\la+s)$ with $\xi>0$ and $s \ge 0$,
the simplest case of (\ref{eq:class2}), then
\[
  \sigma(\la,-\lambda_n)
= \frac{\la+s}{\la(\la+s-\lambda_n \xi)} 
= \frac{s}{\la(s-\lambda_n \xi)} 
   - \frac{\lambda_n\xi}{(\la-\lambda_n\xi)(\la+s-\lambda_n \xi)},
\]
and hence
\begin{equation}\label{eq:specialcn}
c_n(t) = \frac{s}{s-\lambda_n  \xi} 
       - \frac{\lambda_n\xi}{s-\lambda_n\xi}\exp(\lambda_n \xi-s) t.
\end{equation}

\subsection{Admissibility}

In \cite[Theorem~4.3]{jp07} the
following result was established
for the case of a one-dimensional observation, e.g. a point
evaluation (i.e., letting $Y=\CC$).

\begin{theorem}\label{thm:JP-admiss}
Let $a=1+1*k$ with $k\in W^{1,2}(0,\infty)$. Then $C$ is a finite-time
admissible observation operator for (\ref{eqn:uncontrolled})  if and 
only if there are constants $M>0$ and $\omega\in\mathbb R$ such that
\[    
    \sum_{n=1}^\infty \frac{|C\phi_n|^2}
                         {|\lambda|^2|1- \widehat a(\lambda) \lambda_n|^2} 
    \le \frac{M}{\Re \lambda-\omega }, \qquad  \Re \lambda>\omega.
\]  
A similar statement does not hold for infinite-time admissibility, 
see \cite[Example 5.1]{jp07}. 
\end{theorem}

One may rewrite the theorem by duality for the controlled systems
under consideration; however, the kernels given in
Example~\ref{ex:important-kernels} do not satisfy the requirements of
the above result. This observation is a primal motivation for the
present article.

\medskip 

For the control of distributed parameter systems, the study of
one-dimensional inputs may seem a severe restriction of generality.
However, as explained in \cite[Rem.~2.4]{HansenWeiss} for the Cauchy
case, in the case of (finite) $n$-dimensional input spaces,
admissibility is equivalent to the simultaneous admissibility of $n$
one-dimensional systems.  Moreover, the following proposition shows
that a one-dimensional criterion leads immediately to a sufficient
criterion for admissibility of control operators $B: U \to X$ for
infinite-dimensional input spaces. This observation is of great
practical value since the sufficient condition of admissibility in
Theorem~\ref{thm:admissibility-result} can be verified rather easily.
Similar results for the Cauchy problem (i.e., $a\equiv 1$) are well
known in the literature, see e.g. \cite{HansenWeiss, HansenWeiss2}, or
\cite{Haak:Carleson} for the case that $X=\ell_q$ and $U=\ell_p$. Our
proposition generalises directly \cite[Proposition
5.3.7]{TucsnakWeiss}.

\smallskip

Let $U=X=\ell_2$. Let $B: U \to X_{-1}$ be
linear and bounded. Then there are functionals $\varphi_n \in
(\ell_2)^\ast = \ell_2$  such that $(B u)_n = \langle u,\varphi_n 
\rangle$.

\begin{proposition}\label{prop:extension-to-infinite-dim-U}
  Let $X = U = \ell_2$. Let $(\varphi_n)$ be a sequence of elements in
  $U^\ast$ and consider the scalar sequence $b$ defined by $b_n =
  \norm{\varphi_n}$ and let the operator $B$ be defined by $(B u)_n =
  \langle u, \varphi_n\rangle$.

  If $b \in {\mathcal D}(A^\ast)^\ast$ is an admissible input element
  for the resolvent family $(S(t))_{t\ge 0}$, then $B$ is bounded from $U$ to
  ${\mathcal D}(A^\ast)^\ast$ and $B$ is an admissible control
  operator for the resolvent family $(S(t))_{t\ge 0}$ as well.
\end{proposition}
\begin{proof}
  The elementary estimate $|\langle u, \varphi_n
  \rangle| \le \norm{u} \, \norm{\varphi_n}$ implies that $B$ is
  linear and bounded from $U$ to ${\mathcal D}(A^\ast)^\ast$. Now let 
  $u \in L^2(0, \infty; U) = L^2(0, \infty; \ell_2)$ have compact
  support and let $u_j(\cdot)$ denote its coordinate functions. Then
\begin{align*}
    \biggnorm{ \int_0^t S(t{-}s) B u(s)\,ds }_{\ell_2}^2
= & \;\sum_{n=1}^\infty \biggl| \int_0^t c_n(t{-}s) \langle
     u(s), \varphi_n \rangle \,ds \biggr|^2 \\
= & \;\sum_{n=1}^\infty \biggl| \bigl\langle 
      \int_0^t c_n(t{-}s) u(s) \,ds,  \varphi_n\bigr\rangle  \biggr|^2 \\
\le& \;\sum_{n=1}^\infty \norm{\varphi_n}^2 \biggnorm{ 
      \int_0^t c_n(t{-}s) u(s) \,ds}_U^2 \\
 = & \;\sum_{n=1}^\infty \norm{\varphi_n}^2 \sum_{j=1}^\infty
    \biggl| \int_0^t c_n(t{-}s) u_j(s) \,ds \biggr|^2  \\
 = & \;\sum_{j=1}^\infty \biggl[ \sum_{n=1}^\infty  \biggl| \int_0^t
     c_n(t{-}s) \norm{\varphi_n} u_j(s) \,ds \biggr|^2 \biggr] \\
 = & \;\sum_{j=1}^\infty \biggl\|  \int_0^t
     S(t{-}s) b u_j(s) \,ds \biggr\|^2.
\end{align*}
By assumption, $b$ is an admissible input element, and so
\[
\biggnorm{ \int_0^t S(t{-}s) B u(s)\,ds }_{\ell_2}^2
\le \; C \,\sum_{j=1}^\infty \bignorm{ u_j }_{L^2}^2
     = C \,\norm{ u }_{L^2(0,\infty; \ell_2)}^2.
\]
\end{proof}

\bigskip

After this consideration we focus on the case that $U=\CC$, i.e., that
the input space is one-dimensional and let $B\in{\cal L}(U,D(A^*)^*)$.
Then we may write $B=\sum_{n \in \NN} b_n \phi_n$, where
$\{b_n|\lambda_n|^{-1}\} \in \ell_2$.
Further, we assume that the solution family is exponentially stable.\\

We have for $\lambda\in\CC$ with Re$\,\lambda>0$
\begin{equation} \label{eq:action-on-exponentials}
\begin{split}
  {\cal B}_\infty (e^{-\lambda \cdot}) 
  &= \int_0^\infty S(t)Be^{-\lambda t} \, dt \\
  &= \sum_{n=1}^\infty b_n \left( \int_0^\infty c_n(t) 
        e^{-\lambda t}dt \right) \phi_n\\
  &= \sum_{n=1}^\infty b_n \sigma(\la,-\lambda_n)\phi_n \\
  &= \sum_{n=1}^\infty b_n \frac{1}{\lambda-\lambda 
        \widehat a(\lambda) \lambda_n}\phi_n.
\end{split}
\end{equation}
This implies
\[
  \bignorm{ {\cal B}_\infty (e^{-\lambda \cdot}) }^2_{H}
= \sum_{n=1}^\infty \frac{|b_n|^2}{|\lambda|^2 
    | 1-\widehat a(\lambda) \lambda_n |^2}. 
\]
Denote by $k_\la$ the reproducing kernel 
$k_\la(z) := \tfrac{1}{z+\overline{\la}}$. We arrive at the 
following result.

\begin{proposition}
The following is a necessary condition for admissibility of a 
rank-one control operator $B$: There exists a constant $M>0$ such that
\[     
   \sum_{n=1}^\infty \frac{|b_n|^2} { |\lambda|^2 
          |1-\widehat a(\lambda) \lambda_n|^2} 
\le \frac{M}{\Re \lambda},\qquad  \Re \lambda>0,
\]
or, letting $\nu := \sum_n |b_n|^2 \delta_{-\overline{\la_n}}$ where
$\delta_z$ denotes the Dirac mass in $z$, equivalently
\begin{equation}\label{eq:admiss-necessary}
    \norm{ k_{ \frac{1}{\widehat a(\lambda)}} }^2_{L^2(\CC_+, \nu)}
\le M\frac{|\lambda \widehat a(\lambda)|^2}{\Re \lambda}  
    ,\qquad  \Re \lambda>0.
\end{equation}
\end{proposition}

\begin{example}
 \begin{enumerate}
 \item 
  For the particular choice of 
  $\widehat a(\lambda)=\lambda^{-\beta}$, $\beta\in(0,2)$, the necessary 
  condition of the proposition reads
  \begin{equation}\label{eq:klab}
      \bignorm{ k_{\lambda^\beta} }^2_{L^2(\CC_+,\nu)} 
  \le M\frac{|\lambda|^{2-2\beta}}{\Re \lambda},\qquad  \Re \lambda>0.  
  \end{equation}   
 \item In case $\widehat a(\lambda) = \int_0^\infty \tfrac1{\la+s} \,d\al(s)$, the above 
       necessary condition reads
   \[
      \bignorm{ k_{\frac{1}{\widehat a(\la)}} }^2_{L^2(\CC_+,\nu)} 
        \le M  \frac{1}{\Re \la} \biggl|\int_0^\infty \frac{\la}{\la+s} \,d\al(s)\biggr|^2
   \]
 \item If  $\widehat a(\la) = \tfrac1{\log \la}$, the necessary condition reads
   \[
      \bignorm{ k_{\log(\la)} }^2_{L^2(\CC_+,\nu)} 
   \le M  \frac{|\la|^2}{\Re \la |\log \la|^2}
   \]
 \end{enumerate}
\end{example}

There is a strong link between admissibility and Carleson measures in
the Cauchy case $a(t) \equiv 1$, as first observed in \cite{HR83}. For
Volterra systems, we shall establish a similar connection.

\begin{definition}\label{def:embedding--Carleson}
  For $\gamma>0$, a measure $\mu$ on $\CC_+$ is an {\em embedding
    $\gamma$--Carleson measure\/} if for one (and hence all) $q\in
  (1,\infty)$ satisfying $\gamma q>1$ there is an absolute constant
  $M_q$ such that $\norm{f}_{L^q(\CC_+, \mu)} \le M_q
  \norm{f}_{H^{\gamma q}(\CC)}$ for all $f\in H^{\gamma q}(\CC)$.
\end{definition}

We begin with the case $\gamma \in (0,1]$. Here, a measure $\mu$ on
$\CC$ is $\gamma$--Carleson if, and only if there is an absolute
constant $C$ such that
\begin{equation}  \label{eq:geometric--Carleson}
\mu(Q_h)^\gamma \le C \, h  
\end{equation}
for every Carleson square $Q_h$ of side $h$. In case $\gamma=1$ this
characterisation is a celebrated result of Carleson
\cite{Carleson:interpolation, Carleson:corona}, and the extension to
$\gamma<1$ is due to Duren \cite{Duren:Carleson}.  Measures $\mu$ that
satisfy (\ref{eq:geometric--Carleson}) are called {\em geometric
  $\gamma$--Carleson measures}. In the case $\gamma \le 1$ in which
the embedding and geometric $\gamma$--Carleson coincides we simply
speak of $\gamma$--Carleson measures.

\begin{remark}\label{rem:repeated}
  We shall require several times the following easy calculation, where
  we set $\Re\lambda = \xi>0$, and make the substitution $y=\xi t$:
\begin{eqnarray*}
    \norm{ k_\lambda }^p_{H^p} 
&=& \int_{-\infty}^\infty \frac{dy}{(y^2+\xi^2)^{\sfrac{p}2}} \\
&=& \int_{-\infty}^\infty \frac{\xi \, dt}{\xi^p(t^2+1)^{\sfrac{p}2}}
 =   C_p^p \; \xi^{1-p}, \quad \hbox{say}.
\end{eqnarray*}
That is, $\|k_\lambda\|_{H^p} = C_p (\Re \lambda)^{-1/p'}$, where
$C_p$ is a constant depending only on $p$, and $p'$ is the conjugate
index to $p$.
\end{remark}

It is possible to use reproducing kernels as
test functions for the geometric $\gamma$--Carleson property.

\begin{lemma}\label{lem:kernel-positive-lambda}
Assume that for $p,q \in (1,\infty)$ there exists a constant $M>0$
such that
\begin{equation}   \label{eq:kern-abschaetzung}
    \norm{ k_z }_{L^q(\RR^2_+,\mu)} \le M \norm{ k_z }_{H^p(\RR^2_+)}
\end{equation}
for all $z\in \CC_+$. Then $\mu$ is geometric $\gamma$--Carleson for
$\gamma = \sfrac{p}q$. If the support of $\mu$ is contained in a
sector $\cS(\theta)$ with $\theta < \sfrac\pi2$, the conclusion is true
when (\ref{eq:kern-abschaetzung}) merely holds for all $z>0$.
\end{lemma}
\begin{proof}
The proof is a modification of standard arguments that can be seen,
for example in \cite[Lec.~VII]{nik}.

By Remark~\ref{rem:repeated}, $\norm{ k_\la }_{H^p} = C_p
(\Re(\la))^{-\sfrac1{p'}}$. For $\om \in \RR$ and $r>0$, let 
$\la = i\om + r$. Consider the Carleson square $Q_{\om, r}$ with
centre $\om$ and length $r$.
Then the triangle inequality yields $|k_\la(z)| \ge 1/r$ for all $z
\in Q_{\om, r}$ and therefore, 
\begin{align*}
  \mu( Q_{\om, r} ) 
&=   \int_{ Q_{\om, r} } \,d\mu 
\le r^q \int_{ Q_{\om, r} } \bigl| k_\la(z)\bigr|^q \,d\mu \\
&\le M r^q \bignorm{ k_\la(z) }_{H^p}^q 
 =  M C_p r^{q-\sfrac{q}{p'}} = C \; r^{\sfrac1\al}
\end{align*}
This shows that $\mu$ is geometric $\gamma$--Carleson. If $\mu$ has
support in a sector, and if $Q_{\om,h}$ is a Carleson square that
intersects the sector $\cS(\theta)$, then for $x+iy \in Q_{\om,h} \cap
\cS(\theta)$ we have $0 \le x \le h$ and $|y| \le h\tan \theta$. Thus
$Q_{\om,h} \cap \cS(\theta) \subset Q_{0, h \sec\theta}$ which
justifies testing with kernels on the real line. 
\end{proof}

\medskip

\begin{remark}\label{rem:carleson-measure-greater-one}
  The lemma asserts in particular that for $0 < \gamma \le 1$, the
  'reproducing kernel thesis' holds, i.e., if the estimate from
  Definition~\ref{def:embedding--Carleson} holds for the reproducing
  kernel functions $k_\la$, $\la\in \CC_+$, the measure $\mu$ is
  embedding $\gamma$--Carleson. This implication fails in the case
  $\gamma>1$. Indeed, in this case the conditions for a regular Borel
  measure $\mu$ on $\CC_+$ to be embedding $\gamma$--Carleson is
  strictly stronger than to be geometric $\gamma$--Carleson, see e.g.
  \cite{TaylorWilliams} for a concrete counterexample.

  The following necessary and sufficient condition for being embedding
  $\gamma$--Carleson in the case $\gamma >1$ can be found in
  \cite{Videnskii}, see also \cite[Thm.~C]{Lue91}.  Let $S_\mu$ denote
  the \emph{balayage} of $\mu$,
$$
   S_\mu(i \omega) = \int_{\CC_+} p_z(i \omega) d\mu(z),
$$
  where  
\[
   p_z(i\omega)={\pi}^{-1}\frac{\ds x}{\ds x^2+(y-\omega)^2}
\]
  denotes the Poisson kernel for $z=x+iy$ on $i \RR$.  Then
  $\mu$ is embedding $\gamma$--Carleson if
  and only if $S_\mu \in L^{\gamma'}(i \RR)$ where $\gamma'$ is the
  conjugate exponent to $\gamma$.

  A similar characterisation is possible via the Fefferman--Stein
  maximal function $\psi_\mu = \sup_{x\in Q} \frac1h \mu(Q)$
  associated with the measure $\mu$, see \cite{Videnskii}.  The
  arguments of the preceding lemma show that if the measure $\mu$ is
  geometric $\beta$--Carleson and supported in a sector $\cS(\theta)$
  with $\theta<\tfrac{\pi}2$, the Carleson square length $h \ge c |x|$
  and so $\psi_\mu \in L^{\beta', \infty}$ where $\beta'$ is the
  conjugate exponent of $\beta$. Consequently, one obtains from the
  Marcinkiewicz interpolation theorem that  if $\mu$ satisfies
  (\ref{eq:geometric--Carleson}) with $\gamma$ equals $\beta_1$ and $\beta_2$
  then $\mu$ is embedding $\gamma$--Carleson for all $\gamma \in
  (\beta_1, \beta_2)$ (see also \cite{Haak:Carleson}).
\end{remark}

We now introduce the machinery of frames in order to analyse admissibility.

\begin{definition}\label{def:frame}
  Let $H$ be a separable Hilbert space and suppose a sequence
  $(f_n)_{n\ge 1}$ is given.  Then $(f_n)_{n\ge 1}$ is called a {\em
    frame} if there exist constants $B>A>0$ such that
\[
     A \,\bignorm{\varphi}_H^2 
 \le \sum_{n=1}^\infty \bigl|\langle \varphi, f_n\rangle_H \bigr|^2
 \le B \,\bignorm{\varphi}_H^2
\]
for all $\varphi\in H$.
\end{definition}

We recall some basic facts from \cite[Chapter 3]{Daubechies}.  If
$(f_n)_{n\ge 1}$ is a frame, then the so-called {\em frame operator}
$F: H \to \ell_2$, given by $(F \varphi)_n = \langle \varphi, f_n
\rangle$ is clearly bounded. From the very definition of $F$ it
follows that $F^* F$ is bounded and invertible and it can be shown the
elements $\widetilde f_n = (F^* F)^{-1} f_n$ form another (so-called
{\em dual}) frame satisfying
\[
     B^{-1} \,\bignorm{\varphi}_H^2 
 \le \sum_{n=1}^\infty \bigl|\langle \varphi, \widetilde f_n\rangle_H\bigr|^2  
 \le A^ {-1} \,\bignorm{\varphi}_H^2,
\]
together with $\varphi = \sum_n f_n \langle \varphi, \widetilde f_n
\rangle$ for $\varphi \in H$ (see e.g. \cite[Proposition
3.2.3]{Daubechies}). In particular, we may always find a decomposition
$\varphi = \sum c_n f_n$ satisfying the `Besselian' estimate
\[
\bignorm{ (c_n) }_{\ell_2}^2 \le A^{-1} \biggnorm{ \sum_{n=1}^\infty c_n f_n }_H^2 
\]
This elementary property will be used in the proof below. Recall (see
\cite[Def.~I.3.2]{pr}) that the kernel $a$ is {\em sectorial of angle}
$\theta < \sfrac\pi2$ if $|\arg \widehat{a}(\lambda)| \le \theta$ for all
$\lambda$ with $\Re \lambda > 0$, and that the kernel $a$ is {\em
  1-regular} if there is a constant $c>0$ such that $|\lambda
\widehat{a}'(\lambda)|\le c|\widehat{a}(\lambda)|$ for all $\lambda$ with $\Re
\lambda > 0$ (see \cite[Def.~I.3.3]{pr}).

\begin{theorem}\label{thm:admissibility-result}
For $U=\CC$ consider the control operator $B\in{\cal L}(U,D(A^*)^*)$ and the
associated measure  $\mu := \sum |b_n|^2 \delta_{-\lambda_n}$ to the
system $(A,B)$. 
\begin{enumerate}
\item \label{item:B-necess}
 Suppose that $-\lambda_n \in \cS(\theta)$ for all $n\in \NN$ and some 
 $\theta<\sfrac\pi2$, that the kernel $a$ satisfies 
 $\widehat{a}((0,\infty)) = (0,\infty)$ and  
 $|\widehat a(\lambda)| \le C\,|\lambda|^{-\beta}$ for some $C,\beta>0$
 and every $\lambda >0$.

 Then $\mu$ being geometric $\beta$--Carleson is necessary for
 admissibility of $B$.
\item \label{item:B-suff} Suppose that the kernel $a$ is $1$-regular,
  sectorial of angle $\theta <\sfrac\pi2$ and that $|\widehat
  a(\lambda)| \ge c\,|\lambda|^{-\beta}$ for some constants $c>0$ and
  $\beta>\sfrac12$ and every $\la>0$.  Let
  $\beta_2>\beta>\beta_1>\max\{\onehalf, \sfrac\beta3\}$.
 
 Then  $\mu$ being embedding $\beta_1$ and $\beta_2$--Carleson is sufficient for  admissibility of $B$.
\end{enumerate}
\end{theorem}
\begin{proof}
\ref{item:B-necess} Let $B$ be an admissible control operator.  It
follows from Remark~\ref{rem:repeated} on letting
$\sfrac1p+\sfrac1{p'}=1$, that we have
\begin{equation}\label{eq:Hp-norm-of-kernel}
    \norm{ k_{\frac{1}{\widehat a(\lambda)} } }_{H^p} 
 = C_p \bigl(\Re \tfrac1{\widehat a(\lambda)}\bigr)^{-\sfrac1{p'}}
\end{equation}
and hence we obtain
$ \norm{ k_{\frac{1}{\widehat a(\lambda)} } }_{H^p} 
  \ge  C_p \; |\widehat a(\lambda)|^{\sfrac1{p'}}$
for $\Re \lambda>0$. Using condition (\ref{eq:admiss-necessary})  we
have for $\lambda > 0$ 
\[
      \norm{ k_{\frac{1}{\widehat a(\lambda)}}}_{L^2(\CC_+, \nu)} 
\le \sqrt{M} \; \frac{ |\la \widehat a(\la) |}{ \la^\sfrac12 }\\
\le \sqrt{M} C^{1/p}\; |\lambda|^{\sfrac12-\sfrac{\beta}{p}} \, |\widehat a(\la)|^{\sfrac1{p'}} 
 \le M' \; \norm{ k_{\frac{1}{\widehat a(\lambda)}} }_{H^p},
\]
where $p=2\beta$. It follows from Lemma~\ref{lem:kernel-positive-lambda}
that $\mu$ is a geometric $\beta$--Carleson measure, 
which proves the first assertion.\\

\ref{item:B-suff} We now assume that $\mu$ is embedding $\beta_1$ and
$\beta_2$--Carleson and let $p_1=2\beta_1$, $p_2=2\beta_2$. Moreover, let
$u_\lambda(t)= 2 (\Re\,\lambda)^{3/2} te^{-\lambda t}$ and $\mu_{j,k}:= 2^{-j}+i k
2^{-j}$ for $j,k\in \ZZ$.  Then, $\norm{ u_\lambda}_{L^2(0,\infty)}=1$
for all $\la$ with $\Re \la>0$ and in \cite{dupa98} it is shown that
the system $(u_{\mu_{j,k}})_{j,k\in\ZZ}$ is a frame for
$L^2(0,\infty)$.  An easy calculation shows that
\[ 
   B_\infty(u_\lambda) 
=  2 (\Re\,\lambda)^{3/2} \sum_{n=1}^\infty b_n 
    \frac{1-\widehat{a}(\lambda)\lambda_n -\lambda \widehat{a}'(\lambda)\lambda_n}
         {(\lambda -\lambda \widehat{a}(\lambda)\lambda_n)^2} \phi_n.
\]
We further define
\[ 
  g_\lambda(s) := 2 (\Re\,\lambda)^{3/2} \; \frac{1+\widehat{a}(\lambda)s 
                +\lambda \widehat{a}'(\lambda)s}{(\lambda +\lambda \widehat{a}(\lambda)s)^2}.
\]
Using the $1$-regularity of the kernel $a$, there exists a constant $C_1>0$ such that
\begin{eqnarray*}
    \left| \frac{1+\widehat{a}(\lambda)s +\lambda \widehat{a}'(\lambda)s}{(1+ \widehat{a}(\lambda)s)^2}\right|
&=& \left| \frac{\lambda\widehat{a}'(\lambda)+\widehat{a}(\lambda)}{\widehat{a}(\lambda)}\, 
           \frac{1}{1+\widehat{a}(\lambda)s} - \frac{\lambda \widehat{a}'(\lambda)}{\widehat{a}(\lambda)}\,
           \frac{1}{(1+\widehat{a}(\lambda)s)^2}\right|\\
&\le& C_1 \max\left\{ \frac{1}{|1+\widehat{a}(\lambda)s|},  \frac{1}{|1+\widehat{a}(\lambda)s|^2}\right\}.
\end{eqnarray*}
Moreover, $1$-regularity of the kernel implies 
$|\widehat a(\la) | \sim |\widehat a(|\la|) |$ up to a positive constant, see \cite[Lemma 8.1]{pr}.   
Thus $|\widehat a(\lambda)| \ge c'\,|\lambda|^{-\beta}$ 
for all $\la \in \CC_+$. This estimate, together with sectoriality of
the kernel implies for $p\in\{p_1,p_2\}$
\begin{align*}
      \lefteqn{\norm{ g_{\lambda} }_{H^p(\CC_+)}}\\
&\le C_1 \frac{(\Re\,\lambda)^{3/2}}{|\lambda|^2} \max\left\{ 
          \Bignorm{ \frac{1}{1+\widehat{a}(\lambda)s} }_{H^p(\CC_+)}, 
          \Bignorm{ \frac{1}{(1+\widehat{a}(\lambda)s)^2}}_{H^p(\CC_+)}\right\}\\
&\overset{(\ref{eq:Hp-norm-of-kernel})}{=}
      C_2 \frac{(\Re\,\lambda)^{3/2}}{|\lambda|^2} \max\left\{ \frac{1}{|\widehat{a}(\lambda)|}
          \left(\!\Re\, \frac{1}{\widehat{a}(\lambda)}\right)^{-1+1/p}, 
          \frac{1}{|\widehat{a}(\lambda)|^2}\left(\Re\, \frac{1}{\widehat{a}(\lambda)}\right)^{-2+1/p}\right\}\\
& =  C_2 \frac{(\Re\,\lambda)^{3/2}}{|\lambda|^2} \max\left\{ 
          \frac{(\Re\,\widehat{a}(\lambda))^{-1+1/p}}{|\widehat{a}(\lambda)|^{-1+2/p}}, 
          \frac{(\Re\,\widehat{a}(\lambda))^{-2+1/p}}{|\widehat{a}(\lambda)|^{-2+2/p}}\right\}\\
&\le C_3 \frac{(\Re\,\lambda)^{3/2}}{|\lambda|^2} |\widehat{a}(\lambda)|^{-1/p}
\; \le \quad 
     C_4  \frac{(\Re\, \lambda)^{3/2}}{|\lambda|^2} |\lambda|^{\beta/p},
\end{align*} 
for positive constants $C_1$, $C_2$, $C_3$ and $C_4$. 
Thus, 
\[  
   \bignorm{ g_{\mu_{j,k}} }_{H^p(\CC_+)}^2 
\le C_4^2 \, \frac{ 2^{j(1-2\beta/p)} }{(1+k^2)^{2-\beta/p}},
\]
and 
\[ 
M_1 := \biggl(\sum_{j\ge 0,k\in\ZZ} \bignorm{  g_{\mu_{j,k}}}^2_{H^{p_1}}\biggr)^{\onehalf} +
     \biggl(\sum_{j< 0,k\in\ZZ} \bignorm{  g_{\mu_{j,k}}}^2_{H^{p_2}}\biggr)^{\onehalf} < \infty.
\]
Let $u$ be a finite linear combination of the functions $(u_{\mu_{j,k}})_{j,k\in\ZZ}$ and
let  $\alpha_{j,k} = \langle u, \widetilde u_{j,k} \rangle$. By the Besselian property of
the coefficients $\alpha_{j,k}$ we have
\begin{equation}\label{eq:bessel-estimate}
   \sum_{j,k} |\alpha_{j,k}|^2 \le M_2 \,\bignorm{ u}^2_{L^2(0,\infty)}
\end{equation}
for some constant $M_2 > 0$, independent of $u$. This implies
\begin{align*} 
&    \norm{  B_\infty(u) }_H \\
& =  \biggl(\sum_n |b_n|^2 \biggl|\sum_{j,k\in\ZZ} \alpha_{j,k} g_{\mu_{j,k}}(-\lambda_n)\biggr|^2 \biggr)^{\onehalf}\\
& =  \biggl( \int_{\CC_+}   \biggl|\sum_{j,k\in\ZZ} \alpha_{j,k} g_{\mu_{j,k}}(s)\biggr|^2 d\mu(s) \biggr)^{\onehalf}\\
&\le \biggl( \int_{\CC_+}   \biggl|\sum_{j\ge 0,k\in\ZZ} \alpha_{j,k} g_{\mu_{j,k}}(s)\biggr|^2 d\mu(s) \biggr)^{\onehalf}
    +\biggl( \int_{\CC_+}   \biggl|\sum_{j<0,k\in\ZZ} \alpha_{j,k} g_{\mu_{j,k}}(s)\biggr|^2 d\mu(s) \biggr)^{\onehalf}\\
&\le M_3 \biggl(\biggnorm{ \sum_{j\ge 0,k\in\ZZ} \alpha_{j,k} g_{\mu_{j,k}} }_{H^{p_1}(\CC_+)} 
    + \biggnorm{ \sum_{j<0,k\in\ZZ} \alpha_{j,k} g_{\mu_{j,k}} }_{H^{p_2}(\CC_+)} \biggr)
\end{align*}
{}since $\mu$ is embedding $\beta_1$ and $\beta_2$--Carleson. Now since $p_1, p_2 \ge 1$, the
Minkowski and Cauchy--Schwarz inequalities together with (\ref{eq:bessel-estimate}) yield
\begin{align*} 
&\le  M_2 M_3 \;\norm{ u}_{L^2(0,\infty)} \biggl( \biggl(\sum_{j\ge 0,k\in\ZZ} \norm{  g_{\mu_{j,k}}}^2_{H^{p_1}}\biggr)^{\onehalf} +
   \biggl(\sum_{j< 0,k\in\ZZ} \norm{  g_{\mu_{j,k}}}^2_{H^{p_2}}\biggr)^{\onehalf} \biggr)\\
&=  M_1 M_2 M_3 \;\norm{ u}_{L^2(0,\infty)},
\end{align*}
and the proof is done.
\end{proof}

\begin{remark}
\begin{enumerate}
\item If the spectrum of $A$ is contained in a sector
  $\cS(\theta)$ with $\theta<\sfrac\pi2$, the assumption of $\mu$
  being embedding $\beta$--Carleson can be weakened to geometric
  $\beta$--Carleson in the second part of the Theorem by
  using Remark~\ref{rem:carleson-measure-greater-one}.
\item The sectoriality of $a$ with angle $\theta$ already
  implies an estimate on the growth $|\widehat a(\lambda)|$ in the
  half plane. Indeed, as explained in Monniaux--Pr\"uss
  \cite[Proposition 1]{monpr}, 
\[
\log \widehat a(\la) 
= k_0 + \tfrac{i}\pi \int_{-\infty}^\infty \Bigl[\frac{1-ir\la}{\la-ir}\Bigr] 
  \arg(\widehat a(ir)) \frac{dr}{1+r^2}
\]
by a Poisson formula applied to $\arg(\widehat a(\la))$. Here, $k_0$
is a suitable real constant. Considering real $\la>1$ the authors
infer a growth bound $|\widehat a(\la)| \ge c |\la|^{-\al}$ with $\al
= \tfrac{\pi}{2\theta}$.  Combining with \cite[Lemma 8.1]{pr} extends
the estimate to $|\la|>1$.  Inside the unit ball this estimate is not
true in general. Let e.g. $a(t) = t^{-\frac12} e^{-t}$. Then $a$ is a
sectorial kernel of type $\tfrac\pi4$ but $\widehat a(\la) =
\bigl(\tfrac\pi{1+z}\bigr)^\frac12$ attains a finite non-zero limit at
the origin.

\smallskip

Assume that the sectoriality angle satisfies $\theta \le
\tfrac\pi{2\beta}$. Then the growth condition on $\widehat a$ in the
second part of Theorem~\ref{thm:admissibility-result} at infinity is
automatic.  It remains however a non-trivial condition on $\widehat a$
in the origin.
\end{enumerate}
\end{remark}

\subsection{Controllability}\label{sec:controllability}

We are now ready to use the techniques of interpolation to give
conditions for controllability of the Volterra system 
(\ref{eqn:controlled}). Again we assume that the solution 
family is exponentially stable.\\

\begin{lemma}
The following formula for ${\cal B}_\infty$ holds.
\[
{\cal B}_\infty u = \frac{1}{2\pi i}
\sum_{n=1}^\infty b_n \int_{i \RR} \frac{ \widehat u(\la)}{\la_n - 
     (1/\widehat a(-\la))} \frac{1}{\la \widehat a(-\la)} \,d\la \, \phi_n.
\]
\end{lemma}
\begin{proof}
\begin{eqnarray*}
     {\cal B}_\infty u 
&= & \int_0^\infty S(t)Bu(t) \, dt \\
&=&  \sum_{n=1}^\infty b_n \int_0^\infty c_n(t) u(t) \, dt \,\phi_n \\
&=&  \frac{1}{2\pi}\sum_{n=1}^\infty b_n \langle \widehat u, 
       \widehat{\overline{c_n}}\rangle_{H^2(\CC_+)}\, \phi_n,
\end{eqnarray*}
by Plancherel's theorem. Thus, for suitably small $\delta>0$, we have
\begin{eqnarray*}
    {\cal B}_\infty u 
&=& \frac{1}{2\pi} \sum_{n=1}^\infty b_n \int_{-\infty}^\infty 
       \widehat u(\delta+i\omega)\overline{ \widehat {\overline {c_n}}(i\omega+\delta)} 
       \, d\omega \, \phi_n \\
&=& \frac{1}{2\pi} \sum_{n=1}^\infty b_n \int_{-\infty}^\infty 
       \widehat u(\delta+i\omega) \sigma(-i\omega-\delta,-\lambda_n) \, d\omega \, \phi_n\\
&=& \frac{1}{2\pi} \sum_{n=1}^\infty b_n \int_{-\infty}^\infty 
       \frac{ \widehat u(\delta+i\omega)}{-\delta-i\omega +(\delta+ i\omega)\widehat a(-\delta-i\omega)
      \la_n}  \, d\omega \, \phi_n\\
&=& \frac{1}{2\pi i} \sum_{n=1}^\infty b_n \int_{\delta+i \RR} 
       \frac{ \widehat u(\la)}{\la_n - (1/\widehat a(-\la))} 
       \frac{1}{\la \widehat a(-\la)} \,d\la \, \phi_n.
\end{eqnarray*}
\end{proof}

Consider now the kernels defined in (\ref{eq:class2}).
For example, we may take $\widehat a(\lambda)=\xi/(\lambda+s)$, where
$\xi \in \RR$ and $s \ge 0$. We then obtain
\[
  {\cal B}_\infty u = \frac{1}{2\pi i}
\sum_{n=1}^{\infty} b_n \int_{\delta+i \RR} \frac{ \widehat u(\la)(-\la+s)}
{\la(\la_n \xi + (\la-s))}
 \,d\la \, \phi_n ,
\]
which can be calculated using the residue formula as
\begin{equation}\label{eq:binfexplicit} 
{\cal B}_\infty u=\sum_{\substack{n=1 \\
    \Re (s-\lambda_n \xi) > 0}}^\infty b_n \frac{\widehat
  u(s-\lambda_n \xi) \lambda_n \xi}{\lambda_n \xi - s} \, \phi_n.
\end{equation} 
The surjectivity of ${\cal B}_\infty$ reduces to an interpolation
problem of the type analysed in McPhail \cite{mcp} (the case $s=0$ and $\xi=1$
being applied to controllability questions in \cite{jp06}). 

We may use McPhail's theorem as expressed in the half-plane version in
\cite{jp06}. Namely, given $(s_n)$ distinct points in $\CC_+$ and
$(\nu_n)$ non-zero complex numbers, one can find a solution in
$H^2(\CC_+)$ to $F(s_n)=\nu_n x_n$ for every $(x_n) \in \ell_2$, if
and only if $\nu=\sum_{n=1}^\infty \frac{(\Re s_n)^2
|\nu_n|^2}{\varepsilon_n^2} \delta_{s_n}$ is a Carleson measure,
where $\varepsilon_n =  \prod_{k \ne n} \left| 
\frac{s_n-s_k}{s_n + \overline s_k} \right|$.\\
 
The following result therefore generalises part of Theorem 3.1 of \cite{jp06}.

\begin{theorem} In the case $\widehat a(\lambda)=\xi/(\lambda+s)$,
  where $\xi \in \RR$ and $s \ge 0$, exact controllability is
  equivalent to the property that
\[
\sum_{n=1}^\infty \frac{|\Re (s-\lambda_n \xi)|^2|\lambda_n\xi-s|^2}{\varepsilon_n^2|b_n|^2|\lambda_n \xi|^2}\,\delta_{\lambda_n\xi-s}
\]
should be a Carleson measure, where
\[
\varepsilon_n = \prod_{k \ne n} \left| \frac{\xi(\lambda_n-\lambda_k)}{2s-\xi(\lambda_n + \overline \lambda_k)} \right|.
\]
\end{theorem}


Likewise, we may obtain conditions for null controllability in time
$\tau$. The following result reduces to part of Theorem 2.1 of
\cite{jp06} in the case $s=0$, $\xi=1$.

\begin{theorem}
In the case $\widehat a(\lambda)=\xi/(\lambda+s)$, where $\xi \in \RR$
and $s \ge 0$,
null controllability in time $\tau >0$ is equivalent to the property that 
\[
\sum_{n=1}^\infty \frac{|\Re (s-\lambda_n \xi)|^2|\lambda_n\xi-s|^2|c_n(\tau)|^2}{\varepsilon_n^2|b_n|^2|\lambda_n \xi|^2}\,\delta_{\lambda_n\xi-s}
\]
should be a Carleson measure,
where
\[
\varepsilon_n = \prod_{k \ne n} \left| \frac{\xi(\lambda_n-\lambda_k)}{2s-\xi(\lambda_n + \overline \lambda_k)} \right|
\]
and $c_n$ is given by (\ref{eq:specialcn}).
\end{theorem}

\begin{proof}
  This follows on observing that the interpolation problem to be
  solved now has the form ${\cal B}_\infty u= \sum_{n=1}^\infty
  c_n(\tau)x_n \phi_n$ where $(x_n)$ in $\ell_2$ is arbitrary, and
  where ${\cal B}_\infty u$ is given in
  (\ref{eq:binfexplicit}). 
\end{proof}

For higher-order rational functions, the interpolation problems that
arise are more complicated and will repay future investigation.  We
now outline some of the issues involved.  For functions $h$ and $\phi$
we define the weighted composition operator $C_{h,\phi}$ by
\[ 
   (C_{h,\phi} \widehat u)(\la) =h(\la) \widehat u(\phi(\la)).
\] 
We assume that $1/\widehat a(-\cdot)$ maps a piecewise smooth curve
$\Gamma$ bijectively onto $i\mathbb R$.  Let $\phi$ denote the inverse
function of $1/\widehat a(-\cdot)$ mapping $i\mathbb R$ onto
$\Gamma$. 
Assuming that there are no singularities of the integrand below between $i\RR$ and $\Gamma$, so that 
(\ref{eq:just1}) and (\ref{eq:just2}) are equivalent, we
have
\begin{eqnarray}
  {\cal B}_\infty u 
&=&  \tfrac1{2\pi i} 
     \sum_{n=1}^\infty b_n  \int_{i \RR} \frac{ \widehat u(\la)}{\la_n - 
     (1/\widehat a(-\la))} \frac1{\la \widehat a(-\la)} \,d\la\,
     \phi_n \label{eq:just1}\\
 &=&  \tfrac1{2\pi i} 
     \sum_{n=1}^\infty b_n  \int_{\Gamma} \frac{ \widehat u(\la)}{\la_n - 
     (1/\widehat a(-\la))} \frac1{\la \widehat a(-\la)} \,d\la\,
     \phi_n \label{eq:just2} \\ 
&=&  \tfrac1{2\pi i}  \sum_{n=1}^\infty b_n  \int_{i\mathbb R} 
     \frac{ \widehat u(\phi(z))}{\la_n - z} \frac{z \phi'(z)}{\phi(z)} \,dz\, 
     \phi_n\nonumber\\
&=&  \sum_{n=1}^\infty b_n\langle P_{H^2} 
     C_{\frac{z \phi'(z)}{\phi(z)}, \phi} u, k_{\la_n}\rangle \,
     \phi_n.\nonumber
\end{eqnarray}

\begin{example}
  Let $\widehat a(\lambda)=\lambda^{-\onehalf}$, i.e., $\phi(z)=z^2$.
  In this situation we obtain
\begin{eqnarray*}
{\cal B}_\infty u  &=& \tfrac1{\pi i}\sum_{n=1}^\infty b_n  \int_{i\mathbb R} \frac{ \widehat u(z^2)}{\la_n - z}  \,dz\, \phi_n,\\
&=& \sum_{n=1}^\infty 2b_n P_{H^2} (C_\phi u)(\lambda_n) \phi_n. 
\end{eqnarray*}

Now, if $v \in H^2(\CC_+)$ and $v(s)=O(s^{-2})$ as $|s| \to \infty$,
then it follows by an easy estimate of $\int_{-\infty}^\infty
|v((x+iy)^{1/2})|^2 \, dy$ that the function $u:s \mapsto
v(s^{\onehalf})$ lies in $H^2(\CC_+)$, and thus for such functions $v$
we have
\[
B_\infty u= \sum_{n=1}^\infty 2b_n v(\lambda_n) \phi_n.
\]
Thus exact controllability is linked to the condition that for all
sequences $(x_n) \in \ell_2$ there is a function $u \in H^2(\CC_+)$
with $b_n u(\lambda_n^\onehalf) = x_n$, which can once more be
expressed in terms of Carleson measures.
\end{example}

\section{Examples} \label{sec:examples}
The monograph of Pr\"uss~\cite{pr} contains numerous examples of
Volterra systems to which Theorem \ref{thm:admissibility-result} can
be applied. Here we study one particular example.

  Consider a simplified problem of heat conduction with memory in a
  bounded domain $\Omega \subseteq \RR^d$. The uncontrolled situation
  has been studied by Zacher \cite{Zacher}. Integrating his equations
  from zero to $t$ one obtains
\begin{equation}  \label{eq:heat-equation-with-kernel}
     \displaystyle x(t) 
    \displaystyle = 
    \displaystyle x_0 + \int_0^t a(t-s) \Delta x(s)\,ds   + \int_0^t B u(s)\, ds, \qquad t\ge 0
\end{equation}
with some boundary conditions (to be specified later) for the unknown
temperature $x$. Here, the kernel is given by $a(t) = t^\alpha$ where
$\alpha\in [0,1)$ is a material parameter and $B: U \to D(A^*)^*$ is 
the control operator. Notice that the case $\alpha=0$ corresponds to the
classical heat equation and that the (excluded) parameter $\alpha=1$
would correspond to a wave equation.

\smallskip

{\bf Dirichlet boundary conditions.} Consider
(\ref{eq:heat-equation-with-kernel}) under the Dirichlet boundary
condition $x|_{\partial \Omega} = 0$. The problem then reads as
(\ref{eqn:controlled}) where $A$ is the Dirichlet Laplacian.  It is
well known that $A$ is self-adjoint and has compact resolvent and thus
generates a diagonal semigroup. Notice that $\widehat a(\la) =
\Gamma(1{+}\alpha) \la^{-1-\alpha}$ satisfies the growth conditions of
Theorem~\ref{thm:admissibility-result} (here $\Gamma$ denotes the
Gamma function). Moreover, $a$ is evidently $k$--regular for all $k\in
\NN$. From
\[
  H(\la) = (I-\widehat a(\la) A)^{-1} / \la = \tfrac1\la \;
  \la^{\alpha+1} \bigl( \la^{\alpha+1} -  \Gamma(\alpha{+}1) A \bigr)^{-1}
\]
and the sectoriality of $A$ (actually with arbitrary small angle) one
concludes finally that the equation is parabolic in the sense of
\cite[Definition I.3.1]{pr}. Therefore, \cite[Theorem I.3.1]{pr}
assures the existence of a bounded resolvent family $(S(t))_{t\ge 0}$ that is
even $C^\infty( (0,\infty), B(X))$.  \smallskip

We study a rank one control $B: \CC \to D(A^\ast)^\ast$, i.e., $B = (b_n)$ 
via Theorem~\ref{thm:admissibility-result} in the most easy case of a
one-dimensional rod of length one, say $[0,1]$. Then
$\lambda_n=-n^2\pi^2$ for $n\in \NN_{\ge 1}$.
In virtue of the preceding remark it sufficient to verify that $\mu$
is geometric $\beta_1$ and $\beta_2$--Carleson. Let $b_n$ satisfy
$|b_n| \le C n^\delta$ for some $\delta > 0$. Notice that $\lfloor
\sqrt{h}/\pi \rfloor = 0$ when $h\in [0,\pi^2)$. We may thus restrict
ourselves to $h > \pi^2$ in the following estimate:
\[
    \mu( Q_h) 
\le C \sum_{j=1}^{\lfloor \sqrt{h}/\pi \rfloor} j^{2\delta}
\le C \int_{0}^{\lfloor \sqrt{h}/\pi \rfloor} x^{2\delta} \,dx
 =  \frac{C}{1+2\delta} \left(\bigl\lfloor \sqrt{h}/\pi 
     \bigr\rfloor\right)^{1+2\delta}.
\]
Since $h>1$ it is sufficient to establish the estimate
(\ref{eq:geometric--Carleson}) for the maximum of $\beta_1$ and
$\beta_2$ that may be chosen arbitrarily near to
$\beta=1+\alpha$. Using the trivial inequality $\lfloor x \rfloor \le x$, 
one concludes that for $\al \in [0,1)$ given, all elements
$b=(b_n)$ that satisfy
\[
|b_n| \le C n^{\delta} \qquad \text{with}\qquad  
    \delta < \tfrac12 \,\tfrac{1-\al}{1+\al}
\]
are admissible. Let $X=\ell_2$ and let $X_{\theta}$ denote the fractional
domain space of $X$
.  Since $A$ is boundedly invertible, we have $b \in X_{-\theta}$ if, 
and only if
\[
\sum_{n=1}^\infty \left| \frac{b_n}{\la_n^\theta} \right|^2 < \infty
\]
which implies $|b_n| \le C n^{2\theta-\frac12}$. Consequently, all
elements $b \in X_{-\theta}$ with $\theta < \frac12 \frac1{1+\alpha}$
are admissible. This result matches well with the case $\alpha=0$ of
the classical heat equation, where it is well known that $\theta
<\sfrac12$ is sufficient for admissibility  -- a result that fails
 for $\theta=\sfrac12$ (see \cite{Weiss:Conj}).
%

\smallskip

{\bf Neumann boundary conditions.}  Next, we study
(\ref{eq:heat-equation-with-kernel}) under a Neumann boundary
condition $\frac{\partial}{\partial \nu} x = 0$. To this end, let
$\Omega \subset \RR^d$ be a bounded domain that admits a bi-Lipschitz
map from $\Omega$ with constant $L$ onto the unit ball in $\RR^n$.
Let $A$ denote the negative Neumann Laplacian $-\Delta_N$ on $\Omega$. The
eigenvalues of $A$ can be arranged according to their multiplicities
as 
\[
0 =\mu_0 < \mu_1 \le \mu_2 \le \ldots
\]
As in the one-dimensional case with Dirichlet boundary conditions, the
uncontrolled problem is parabolic and the existence of a bounded
solution family $(S(t))_{t\ge 0}$ is assured by \cite[Theorem I.3.1]{pr}. It
is known (see e.g.  \cite{Kroeger}) that the eigenvalues of $\Delta_N$
satisfy an estimate
\begin{equation}  \label{eq:neumann-laplace-eigenvalues}
C_{|\Omega|, d} \bigl(1- C_{d, L} \,n^{-\sfrac1d} \bigr) \, n^\sfrac2d
\le \mu_n
\le C_{d, L}' \, n^\sfrac2d.
\end{equation}
where the subscripts refer to the dependencies of the constant on the
dimension $d$, the bi-Lipschitz constant $L$ and the volume $|\Omega|$
respectively. Let $\phi_{n}$, $n\ge 0$ be a basis of eigenvectors of
$\Delta_N$, normalised in a suitable Hilbert function space $X$, 
and let $B: U \to D(A^\ast)^\ast$ be bounded, where $U$ is another Hilbert
space. Then
\[
B u = \sum_{n=1}^\infty  \phi_{n} \bigl\langle u,f_{n} \bigr\rangle_{U}
\]
with $f_{n} = B^\ast \phi_n \in U^\ast$. Combining
Theorem~\ref{thm:admissibility-result} and
Proposition~\ref{prop:extension-to-infinite-dim-U} yields that $B$
is admissible provided that the measure
\[
\nu = \sum_{n=1}^\infty \delta_{\mu_n} \norm{ f_{n} }_{U}^2 
\]
is geometric $\beta_1$ and $\beta_2$--Carleson for suitable $\beta_i$
close to $\beta=1+\al$. Notice that we did not put any weight at
$\mu_0 = 0$ since this would destroy the $\beta$--Carleson property
for all $\beta>0$.  Since the support of $\nu$ does not intersect with
a small ball around the origin, it is again sufficient to establish
estimate (\ref{eq:geometric--Carleson}) for the maximum of $\beta_1$
and $\beta_2$; moreover, we may restrict to large $\mu_n$, i.e.,  we
may assume without loss of generality that $\tfrac1c n^{\sfrac2d} \le
\mu_n \le c \,n^{\sfrac2d}$ for some $c>0$.  Then, essentially the
same calculation as in the example of the one-dimensional rod above
yields that if $ \norm{ f_{n} }_{U} \le C n^\delta$ for all $n\in
\NN$, $B$ is an admissible control operator provided that $\delta <
\frac12 \frac{ \sfrac2d -1 -\al}{1+\al}$.

In the linear case, multiplying a boundary control system of the form
$x'(t) + \Delta x(t) =0$, $\frac{\partial}{\partial \nu} x = u$ with a
test function and integrating by parts shows that it fits into the
abstract setting $x'(t) + Ax(t) = B u(t)$ where $A=\Delta_N$ and $B$
is the adjoint operator of the Dirichlet trace, see e.g. \cite{BGSW}.
One needs therefore to estimate the boundary traces of the Neumann
eigenvectors. 
The choice of the Hilbert function spaces $X$ (in the domain
$\Omega$) and $U$ (on its boundary) plays an important r\^ole.

Take e.g. $X= H^1(\Omega)$ and $U=L^2(\partial \Omega)$ and assume
that $\Omega$ has a $C^1$ boundary. Together with $\mu_n \ge c
\,n^{\sfrac2d}$ for large $n$ and the boundedness of the Dirichlet
trace from $H^{\sfrac12}(\Omega)$ to $L^2(\partial\Omega)$, it readily
follows in the case $d=1$ and $d=2$ that for $\al \in [0,1)$ the
boundary control $B$ is admissible.  For $d \ge 3$ the problem can be
analysed by passing to a higher-order Sobolev space. For smoother
domains, it is possible to use more sophisticated estimates for the
Dirichlet trace operator, such as those given by Tataru~\cite{tataru}.

\end{document}